\documentclass[12pt,titlepage]{article}
\usepackage{amsmath}
\usepackage{amsfonts}
\usepackage{amssymb}
\newtheorem{satz}{Theorem}

\newtheorem{prop}[satz]{Proposition}
\newtheorem{rem}[satz]{Remark}
\newenvironment{proof}[0]{PROOF.}{\hfill $\Box$ \bigskip \newline}

\def\span{\mbox{\rm span}}
\def\supp{\mbox{\,\rm supp\,}}
\def\y{\lbrack\!\lbrack} 
\def\yy{\rbrack\!\rbrack}

\def\h1h{H^1(H)}

\def\h1x{H^1(X, d, \mu )}

\font\textmsbm= msbm10 scaled 1200
\font\scriptmsbm= msbm7 scaled 1200
\font\scriptscriptmsbm= msbm5 scaled 1200

\font\textmsbm= msbm10 scaled 1200
\font\scriptmsbm= msbm7 scaled 1200
\font\scriptscriptmsbm= msbm5 scaled 1200

\def\BMO{\mathchoice{\mbox{BMO}}{\mbox{BMO}}{\mbox{\scriptsize
BMO}}{\mbox{\tiny BMO}}}
\def\Id{\mathchoice{\mbox{Id}}{\mbox{Id}}{\mbox{\scriptsize
Id}}{\mbox{\tiny Id}}}

\def\proof{\goodbreak\noindent{\sc Proof. }\nobreak}

\def\endproof{\par\nobreak\hbox to \hsize{\hfil\vrule width 5pt height
5pt}\goodbreak\vskip 3pt}

\newfam\tmwfama
\textfont\tmwfama\textmsbm
\scriptfont\tmwfama\scriptmsbm
\scriptscriptfont\tmwfama\scriptscriptmsbm
\def\bbb{\mathbb}

\def\bR{{\bbb R}}

\def\bN{{\bbb N}}

\def\bN{{\bbb N}}

\def\bS{{\bbb S}}

\def\bZ{{\bbb Z}}

\def\s{\sigma}

\def\t{\tau}

\def\sm{\setminus}

      \def\sbe{\subseteq}

\def\la{\langle}
\def\ra{\rangle}

\def\cB{{\cal B}}
\def\cC{{\cal C}}
\def\cD{{\cal D}}
\def\cE{{\cal E}}
\def\cF{{\cal F}}

\def\cH{{\cal H}}

\def\cL{{\cal L}}

\def\Id{\mbox{ \rm Id} }

\title{Extrapolation of Vector valued 
Rearrangement Operators II }
\author{Paul F.X. M\"uller
\thanks{Supported by the Austrian Science Foundation P 20166-N18,
\newline
2000 MSC 46B42, 46B70,47B37.}}
\begin{document}
\maketitle
\begin{abstract}
 We determine the extrapolation law of rearrangement operators
acting on the Haar system in vector valued $H^p$ spaces: 
If  $ 0<q \le p <2 , $ then, 
$$
 \|  T_{\t, q } \otimes \Id_X \| _{q}^{\frac{q}{2-q}}  \le 
A(p,q) \|  T_{\t, p } \otimes \Id_X \| _{p}^{\frac{p}{2-p}} .
$$
For a fixed Banach space $X, $ the extrapolation range $ 0<q \le p <2  $ is 
optimal.
If, however, there exists $ 1 < p_0 < \infty , $ so that 
$$
\|  T_{\t, p_0 } \otimes \Id_E \| _{L^{p_0}_E} < \infty ,\quad\text{for each UMD space E,}
$$
then for any  $1 < p <  \infty $,
$$
\|  T_{\t, p } \otimes \Id_E \| _{L^{p}_E} 
< \infty,
$$
 for any UMD space $E.$ (The value   $p_0 =2$ is {\it not} excluded.)
\end{abstract}

\section{Introduction}
In this  note we identify the extrapolation law of rearrangement operators
acting on the Haar system in vector valued $H^p$ spaces: 
Rearrangement operators are given by an 
injective map $\tau $ acting on dyadic intervals. 
We have 
$$T_{ \t,p} \left( \frac{ h_I}{|I|^{1/p}}\right) =  \frac{h_{\t (I) }}{|\t(I)|^{1/p}}, \quad  0 <  p <2,$$
and by linear extension obtain an operator on the span of the Haar system 
If  $ 0<q \le p <2 , $ then, we prove that   
\begin{equation}\label{23jan}
 \|  T_{\t, q } \otimes \Id_X \| _{q}^{\frac{q}{2-q}}  \le 
A(p,q) \|  T_{\t, p } \otimes \Id_X \| _{p}^{\frac{p}{2-p}} ,
\end{equation}
as well as
\begin{equation}\label{23jan1}
 \|  T_{\t, q } \otimes \Id_X \| _{q}\le A(p,q) \|  T_{\t, p } \|_p^{2\left(\frac{1}{q} - \frac{1}{p}\right)}
\|  T_{\t, p } \otimes \Id_X \| _{p}.
\end{equation}
By arithmetic \eqref{23jan} is implied by \eqref{23jan1}. 
We define in this paper (see Section~\ref{preliminary}) the norm in  vector valued $H^p$ spaces as  the $L^p$ norm of 
Rademacher averages, and use
systematically the notation 
$$  \|  T_{\t, p } \|_p =  \|  T_{\t, p } : H^p \to H^p\| , $$
$$
\|  T_{\t, p } \otimes \Id_X \| _{p} = \|  T_{\t, p } \otimes \Id_X :
H_X^p \to H_X^p \|   $$
and 
$$
\|  T_{\t, p } \otimes \Id_X \| _{L^p_X} = \|  T_{\t, p } \otimes \Id_X :
{L^p_X} \to {L^p_X} \|.   $$
We put the result of this paper,  \eqref{23jan} and  \eqref{23jan1}, into perspective  by reviewing its predecessors.
\paragraph{Scalar valued extrapolation.} 

The extrapolation law \cite{gmp} for scalar valued 
rearrangement operators on dyadic $H^p $ spaces is this, 
\begin{equation}\label{15jan081} 
a_{p,q} \|T_{\t,q}\|_q^{q/( 2-q)}  \le \|T_{\t,p}\|_p^{p/( 2-p)} \le A_{p,q} \|T_{\t,q}\|_q^{q/( 2-q)}, \quad \quad 0 < q \le  p < 2 . 
\end{equation}
Thus, boundedness of $ T_{\t, p }$ on $H_p$  for one value of 
$0< p < 2 $ implies boundedness of   $ T_{\t, q } $  on $H_q$ for all 
values of $q$ where $  0<q  <2.$ 
Separately, the boundedness of $T_{\tau, 1}$ on $H^1$ is equivalent to 
$ \s = \t^{-1} $ respecting Carleson Constants. Write 
$$  \y\cC \yy = \sup_{I \in \cC } \frac{1}{|I|} \sum_{J \sbe I , J \in  \cC}|J|, $$
then \cite{pfxm} 
\begin{equation}\label{3feb1}
\|T_{\tau, 1}\|_1 \sim \sup \frac{\y\s(\cC) \yy }{ \y\cC \yy} ,
\end{equation} 
where the supremum is extended over all $\cC \sbe \t ( \cD) .$

\paragraph{Vector valued  rearrangements. }
By way of example \cite{gmpre} it is easy to see that
 isolating  intrinsic criteria  
characterizing boundedness of 
$   T_{\t,p} \otimes \Id_X ,$
and the search for extrapolation theorems 
represent two different lines of research, both of which are different 
from the scalar valued setting.

For a  rearrangement $ \t_0 $ defined in \cite{gmpre}
 the boundedness of 
$ T_{\t_0 ,p  } \otimes \Id_X ,$ $  1 < p \le  2, $  
implies Rademacher-type $p$ for  $X .$ At the same time the scalar valued extension  of $ T_{ \t_0 ,1}$
is bounded on $H^1 $ while its inverse is unbounded. For 
 $ \t_0 $ and $ 1 < p \le  2$ we have \cite{gmpre},
\begin{equation}\label{15jan082} \|  T_{ \t_0 ,1} \|_1 < \infty ,\quad  
\|  T_{ \t_0 ,1}^{-1} \|_1 = \infty  \quad\text{and}\quad  \text{Type} _p (X) \le   C \|  T_{\t_{0},p } \otimes \Id _X \| _{L^p_X}.
\end{equation}
This example puts restrictions on  possible extrapolation theorems 
for vector valued rearrangement operators. For instance 
the  right hand side estimate in \eqref{15jan081}
is ruled out when \eqref{15jan082} holds. (Just recall
 how the  
 Rademacher-type of a Banach space depends on $p.$)
In \cite{gmpre} we defined  $ \t_1$ such  that  the scalar valued extension
$ T_{ \t_1,1 }$ is an isomorphism on $H^1$ and the boundedness of 
$ T_{\t_{1},1 } \otimes \Id _X $  on $L^p_X$ 
implies the UMD condition for $X.$

We formulate now three general extrapolation 
estimates that are not yet ruled out by the examples discussed above.
\begin{enumerate}
\item 
The first concerns the extrapolation of isomorphisms
across the entire scale of vector valued $L^p ( 1 < p < \infty)$ spaces. Let $X$ satisfy the UMD property, assume  $ T_{\tau, 1} $ is a 
scalar valued isomorphism on $H^1$ and 
$$ \|  T_{\t, 2 } \otimes \Id_X \| _{L^2_X} \|  T_{\t, 2 }^{-1} \otimes \Id_X \| _{L^2_X}
< \infty \quad, $$
then 
$$ \|  T_{\t, p } \otimes \Id_X \| _{L^p_X} \|  T_{\t, p }^{-1} \otimes \Id_X \| _{L^p_X}
< \infty , 
\quad 1 < p < \infty .$$
This result is  in \cite{gmpre}  
were the proof is based on geometric and combinatorial characterizations 
of  rearrangements $\t $ when  $ T_{\tau,1} $ is a 
scalar valued isomorphism on $H^1.$  
\item 
The second extrapolation estimate  that is not ruled 
out by the examples discussed above 
asserts the following. If 
$$ \|T_{\t,1}\|_{1} < \infty
\quad\text{and} \quad  \|  T_{ \t,2 } \otimes \Id_X \| _{2} < \infty , 
$$ then 
$$ \|  T_{p,\t } \otimes \Id_X \| _{p} < \infty ,\quad 0 < p \le  2.$$
Related are the estimates
\begin{equation} \label{2812so1}
 \|  T_{ \t,q} \otimes \Id_X \| _{q}^{q/(2-q)} \le A( p,q)
\|  T_{\t, p } \otimes \Id_X \| _{p}^{p/(2-p)} ,\quad 0 < q < p < 2 . 
\end{equation}
\item The above extrapolation estimates were stated for one  fixed Banach space $X.$
The following assumes boundedness of the rearrangement operator for {\em each }
Banach space with the UMD property. Clearly this is a more restrictive hypothesis
on the underlying rearrangement,
so  the resulting extrapolation estimates should be  stronger. 
If there exists $1 < p_0 <\infty $  so that 
\begin{equation} \label{8feb13} 
\|  T_{\t, p_0 } \otimes \Id_E \| _{L^{p_0}_E} < \infty ,\quad\text{for each UMD space E,}
\end{equation}
then for all  $1 < p <  \infty $
\begin{equation} \label{8feb14} 
\|  T_{\t, p } \otimes \Id_E \| _{L^{p}_E} \|  T_{\t, p }^{-1} \otimes \Id_E \| _{L^{p}_E} < \infty,\end{equation}
for each UMD space E. 
(Note that now   $p_0 =2$ is permissible  in the hypothesis \eqref{8feb13}.) 
\end{enumerate}
In this paper  we  prove that \eqref{8feb13} implies \eqref{8feb14}, provide a proof of \eqref{2812so1},  
and of the implication stated before \eqref{2812so1}.
We point out two direct predecessors  to the present  work.
In \cite{gmpre}, we 
applied  Maurey's \cite{ma2}
 extrapolation-by-factorization method 
to  $\t -$monotone operators.  
For 
$
|I|= |\t(I)| $ and  UMD spaces $X$ we proved  
that 
$$ \|  T_{\t, q} \otimes \Id_X \|^{q/(2-q)}_{L^q_X} \le A(p,q,X) 
\|  T_{\t, p } \otimes \Id_X \| _{L^p_X}^{p/(2-p)} ,\quad 1 < q\le  p \le  2 . $$
The second predecessor  is K. Smela's \cite{smela} very recent proof 
of the scalar extrapolation theorems in \cite{gmp}. 
The results of the present paper 
were obtained by comparing the integral estimates 
for the maximal functions 
\begin{equation} \label{2812so1a}
\mu _{\cH } ( t ) = \sup _{I \in \cH } \frac{ |\s (I) |} {|I|} 1_{I }(t) 
\quad\quad
\text{where} \quad \s = \t^{-1} ,\end{equation}
used by K. Smela  \cite{smela}, to  the methods employed 
in \cite{gmp, gmpre}  and \cite{pfxm}. 
\paragraph{Acknowledgement:}
I would very much like to thank Stefan Geiss
for  helpful discussions during the preparation of this work and
for suggestions improving its presentation.
  
\section{Preliminaries} \label{preliminary}
Here we collect frequently used facts and theorems. 
We routinely use  \cite{m} as  reference.
\paragraph{Collections of dyadic Intervals.}
Let $\cD $ denote the collection of (half-open) dyadic intervals
contained in the unit interval
$$
[(k-1)2^{-n} , k 2^{-n}[ , \quad 1\le k \le 2^{n} , \quad n \in \bN . $$
For $n \in \bN $ write  
$\cD_n = \{ I \in \cD : |I| = 2^{-n} \}. $  
For a collection of dyadic intervals $\cE $ we us the $*-$ notation to 
to denote the pointset covered by $\cE $ thus
$$ \cE^* = \bigcup_{I \in \cE } I . $$

Given $I \in \cD $ we denote by $G_1(I, \cE )$ the maximal dyadic intervals that are in $\cE $ and strictly contained in $I .$
Note that by telescoping,  for a dyadic interval $I,$
$$ |I| = \sum_{K \in \cE , K \sbe I} |K| - | G_1(K, \cE)^*| . $$
The $n-$th generation of the dyadic intervals in $\cE$ underneath $K$ is defined inductively as 
$$ G_n (K | \cE ) = \bigcup_{J \in G_{n-1}(K| \cE )} G_1(J | \cE) . $$
Let $\cL $ be a collection of dyadic intervals. We say that $\cB (I)\sbe \cL$
is a block of dyadic intervals in $\cL $ if  the following conditions hold.
\begin{enumerate} 
\item In $\cB (I)$ has a unique maximal interval, namely the interval $I .$
\item If $J \in \cB (I)$ and $ K \in \cL $ then 
        $$ J \sbe K \sbe I \quad\text{implies} \quad K \in  \cB (I) . $$
 \end{enumerate}
\paragraph{The Haar system.}
Denote by $ \{ h_I :  I \in \cD \} $
the $L^\infty-$ normalized Haar system, where $h_I$ is supported on $I$ 
and 
$$
h_I = \begin{cases} 1 \quad &\text{ on the left half of }I ;\\ 
                    -1  \quad &\text{ on the right half of }I . 
      \end{cases}
$$ 
For $ f \in  L^p$ we define its dyadic square function as 
$$ S(f) = \left( \sum_{I \in \cD } \la f,\frac{h_I}{|I|}\ra ^2 1_I \right)^{1/2}.$$
The Marcinkiewicz- Zygmund interpretation of  R.E.A.C. Paley's theorem
asserts that 
$$ c_p \| f\|_{L^p} \le \| S(f)\|_{L^p} \le C_p \| f\|_{L^p} , \quad ( 1< p  <\infty).$$
Given $ 0 < q < \infty $ define dyadic $H^q,$  
to be the completion of $\span \{ h_I :  I \in \cD \} $ 
under the (quasi-) norm 
is given by 
  $$  \| f\|_{H^q} = \| S(f) \|_{L^q} . $$ 
\paragraph{The dual of dyadic $H^1 . $}
Define $ f \in \BMO $ if 
$$
 \|f \|_{\BMO} ^2 = \sup_{I \in \cD }  \frac {1 }{|I|}  
\| \sum_{ J \sbe I }  \la f ,\frac {h_J}{|J|} \ra h_J
\|^2_2 < \infty . $$
Let 
$ \cL = \{ I \in \cD : \la f , h_I \ra \ne 0 \} $ then
$$ \|f\|_2^2 \le \|f\|_{\BMO}^2 |\cL ^*| .$$  
The space $\BMO $ is (identified with) the dual to $H^1 .$ 
The pairing between
$ f \in \BMO $ and $ g \in H^1 $ is 
$$ \la f , g \ra = \lim _{ n \to \infty } \int_0^1 f_n g dt ,
\quad\text{where}\quad 
 f_n = \sum _{\{I : |I| > 2^{-n} \} }  \la f ,\frac {h_I}{|I|}\ra h_I . $$
\paragraph{Kahane's principle of contraction and Kahane's inequality.}
See \cite{kahane},
\cite{w}.
Let $\{r_n \} $ denote the sequence of independent 
$\{+1, -1\}$ valued Rademacher functions. Let 
$x_n \in X $ be a sequence in a Banach space $X$
and let $ a_n \in \bR $ so that $ |a_n| \le 1 . $
Then,  
$$\int_0^1 \| \sum _{n=1}^N  r_n(t) a_n x_n \|_X dt \le 
\int_0^1 \| \sum _{n=1}^N  r_n(t)  x_n \|_X dt . 
$$ 
We apply the  principle of contraction in combination with
the  Kahane's inequality asserting that 
$$
\left(\int_0^1 \| \sum _{n=1}^N  r_n(t)  x_n \|_X ^pdt\right)^{1/p} \le 
C_p
\int_0^1 \| \sum _{n=1}^N  r_n(t)  x_n \|_X dt, \quad  1<p<\infty . 
$$ 
\paragraph{Vector valued dyadic Hardy Spaces .} See \cite{bu1, ms}. \nocite{blasco}
Given a Banach space $X$ and $x_I \in X ,$ 
define $f = ( x_I : I\in \cD )$  to be the $X$ valued vector indexed and 
ordered by the dyadic intervals. 
Define the square function of $f$ 
as
$$ \bS(f)(t) = \lim_{n\to \infty} \left(\int_0^1  \| \sum_{\{I  : |I| \ge 2^{-n}\}} r_I(s) 
x_I  h_I(t) 
\|_X^2 ds\right)^{1/2}, 
$$
where $ \{r_I \} $ is an enumeration of the independent Rademacher system. 
Let $ 0 < p < \infty. $
We say that   $ f \in H^p_X$ if
$$ \|f\|_{H^p_X} = \|\bS(f)\|_{L^p} < \infty. $$
We (should not hesitate to) identify $f = ( x_I : I\in \cD )$ with its formal Haar series
$ f = \sum_{I \in \cD } x_I h_I . $
If  $1 < p < \infty $ and if  $X$ 
 has  the UMD property, 
 $L^p_X$ (the Bochner-Lebesgue space) 
and $H^p_X$ coincide  with equivalent norms.
\paragraph{Rearrangement Operators.} Let  $ 0 < p < \infty. $ 
Assume  $ \t : \cD \to \cD $ is  injective and 
 $f = ( x_I : I\in \cD )$ in  $ H^p_X  .$ 
The collection of dyadic intervals $ \{ I \in \cD : x_I \ne 0 \} $ 
is the Haar support of $f.$
Define the rearrangement operator $T_{\t, p } \otimes \Id_X $ in terms of 
formal Haar series by the relation 
$$ T_{\t, p } \otimes \Id_X(f) = \sum_{I \in \cD } x_I \left(\frac{|I|}{|\t (I)|}\right)^{1/p}h_{\t (I)} . $$
Equivalently, in vector notation,
$$ T_{\t, p } \otimes \Id_X(f) = \left( x_{\s (J)} \left(\frac{|\s (J)|}{| J|}\right)^{1/p} : J \in \t(\cD) \right ), $$
where $ \s = \t^{-1} : \t(\cD) \to \cD .$  We write
$$ \| T_{\t, p } \otimes \Id_X : H^p_X \to H^p_X\| = \sup \left\{
\| T_{\t, p } \otimes \Id_X(f)\|_{H^p_X} \right\}, $$
where the supremum is extended over all 
$f$ in the unit ball of  $ H^p_X$
with finite Haar support. 
   \paragraph{Dyadic Atoms.} Let $ 0 < p < \infty $
and  $b_J \in X . $  Define $a = ( b_J : J \in \cD ) $
to be a dyadic $H^p_X -$ atom
if there exists a dyadic interval $I$ so that 
$$ \supp \bS(a) \sbe I \quad\text{and}\quad  \quad 
\|\bS(a)\|_\infty \le |I|^{-1/p} .$$
Note that 
$ \|a\|_{H^p_X} \le 1 , $
for a  dyadic $H^p_X -$ atom. \paragraph{Atomic decomposition.}
 Let $ 0 < p < \infty $ and fix 
$ f = (  x_J : J \in \cD )$ such  that  $ f \in H^p_X .$ 
We employ the atomic decomposition that results from stopping time arguments
applied to the square function $\bS(f) : $
There  exists a decomposition of $\cD $ into blocks of dyadic intervals 
$\{ \cB(I) : I \in \cE \}$ and integers $n(I)\in \bZ $ so that 
\begin{equation}\label{27jan093}
\sup _{I\in\cE} \frac{1}{|I|} \sum_{J \sbe I , J \in \cE } |J| \le 4 , 
\end{equation}
\begin{equation}\label{27jan092}
\bS(f_I) \le   2^{n(I)}, \quad\text{for}\quad f_I =\sum_{J \in \cB (I) }x_J h_J  ,
\end{equation}
and
\begin{equation}\label{27jan091}
c \|f\|^p_{H^p_X} \le \sum_{I \in \cE } |I| 2^{pn(I)} \le A_p \|f\|^p_{H^p_X}.
\end{equation} 
Define 
$$ \lambda_I = |I|^{1/p} 2^{n(I)} , \quad I \in \cE . $$ 
By \eqref{27jan092}, 
$ a_I = f_I /\lambda_I $
is a dyadic $H^p_X -$ atom so that   
$$ f = \sum_{I \in \cE } \lambda_I a_I  \quad\text{(formal Haar series)} $$
and by \eqref{27jan091}, 
$$     
c_p \|f\|^p_{H^p_X} \le \sum_{I \in \cE } \lambda_I ^p \le A_p \|f\|^p_{H^p_X}.
$$
The atomic decomposition as cited above originates with \cite{JJ}.
The decomposition of $ \cD $ into blocks $ \cB(I) $ is described 
(for instance) in 
\cite[Pages 42-44]{m}; the right hand side estimate of 
\eqref{27jan091} transfers directly to the range $0<p<\infty$ and to 
the square function defining the spaces  $H^p_X . $
For the left hand side estimate of \eqref{27jan091} we distinguish between
the cases $ 0 < p < 1 $ and $ 1 \le   p < \infty . $
For $0 < p  < 1 $ use the quasi-triangle inequality for the spaces 
 $H^p_X . $  For  $ 1 \le  p < \infty  $ exploit \eqref{27jan093}
and adapt the proof of  \cite[Lemma 3]{gmbull} to yield 
$$  \|f\|^p_{H^p_X} \le C \sum_{I \in \cE }  \|f_I\|^p_{H^p_X} \quad 
(1 \le p < \infty),$$
where $C> 0$ depends just on the upper estimate \eqref{27jan093} for the 
Carleson constant of $\cE . $

\section{Extrapolation by Factorization and Carleson measure}
In this section we prove the results of this paper.
 Let $ \t : \cD \to \cD $ be injective with inverse 
$ \s : \t ( \cD ) \to \cD . $
K. Smela   based his proof \cite{smela} on the fact that 
the maximal function $\mu _{\cH } ( t )$
mediates between $ T_{\t, 1} $ and $T_{\t, 2} . $
This is the link between  extrapolation  of  general rearrangement operators 
and  factorization that extends the use of  $\t-$monotone operators 
in \cite{gmpre}. 
We let $ f = \sum _{J \in \cL}  a_J h_J ,$
with $a_J \in \bR . $
For $ t \in [0,1] $ fixed we estimate square functions 
$$
\begin{aligned}
S(  T_{\t, 1} (f) ) (t) &= \left( \sum _{J \in \cL}  a_J^2 1_{\t(J)}(t)
\frac{|J|^2}{|\t(J)|^2} \right)^{1/2} \\
&\le 
\left(  \sup _{J \in \cL } \frac{ |J |} {|\t(J)|} 1_{\t(J)  }(t)  
\right)^{1/2} \left( \sum _{J \in \cL}  a_J^2 1_{\t(J)}(t)
\frac{|J|}{|\t(J)|} \right)^{1/2}. 
\end{aligned}
$$
The right hand side factor coincides with $ S(  T_{\t, 2} (f) ) $
and the left hand side factor may be rewritten as 
$$
\mu _{\cH } ( t ) = \sup _{I \in \cH } \frac{ |\s (I) |} {|I|} 1_{I }(t) ,
\quad\quad
\text{where} \quad  \cH = \t(\cL ) ,\, \s = \t^{-1} . $$
Summing up, for any $f $ with Haar support  $ \cL  $
we have
 the factorization 
$$
S(  T_{\t, 1} (f) ) (t) \le \mu _{\cH } ( t )^{1/2} S(  T_{\t, 2} (f) )(t)
\quad\text{where} \quad   \cH = \t(\cL ) . $$ 
 
To $ \t : \cD \to \cD $  injective with inverse 
$ \s : \t ( \cD ) \to \cD  $ we define $S_\s $ to be the linear extension of 
the map $ S_\s ( h_I)  = h_{\s (I) } $ when $ I \in \t ( \cD )$ and 
 $ S_\s ( h_I)  = 0  $ when $ I \in \cD \sm \t ( \cD ) .$
The content of following proposition 
appeared  in \cite{smela} by K. Smela. We present it here with a
short proof emphasizing  the connections  of  maximal functions
to Carleson Measure and $\BMO . $ Only the right hand side of 
the inequality will be needed later in the extrapolation proof.

\begin{prop}\label{proposition} 
$$
\frac{\| S_\s\|_{\BMO} ^2}{A} \le 
\sup_{ \cH \sbe \t( \cD )} \frac{1}{ |\s(\cH)^*|}
 \int_0^1 
\mu _{\cH } ( t ) dt 
\le \| S_\s\|_{\BMO} ^2 , 
$$
where $A> 0$ depends on the constants of the atomic decomposition for $H^1. $
\end{prop}
\proof 
First estimate the right hand side. 
Fix $ \cH \sbe \cD . $
Without loss of generality we assume 
that $ \cH $ is a finite collection of intervals.
In the first step of the argument we resolve the maximal function.
To $ t \in [ 0, 1 ] $ choose $J_t \in \cH $ so that 
$$  \frac{ |\s (J_t) |} {|J_t|} = \sup _{J \in \cH } \frac{ |\s (J) |} {|J|} 1_{J }(t). 
$$
Put $ \cB = \{ J_t :  t \in [ 0, 1 ] \} \sbe \cH .$
Fubini's theorem yields 
\begin{equation} \label{2812so2}
 \int_0^1 \sup _{J \in \cH }
\frac{ |\s (J) |} {|J|} 1_{J }(t) dt  =
\sum_{K \in \cB }   | \sigma (K) |  \frac{|K|- |G^*_1( K | \cB )| }{|K|} .
\end{equation}
Thus resolving the  maximal functions led us to evaluating Carleson Measure.
Next we obtain  estimates for the right hand side of \eqref{2812so2}.
For $ K \in \cB $ put 
$$ c_K^2 =  
(|K|- |G^*_1( K | \cB )| )/|K| $$ 
and define 
$$ f = \sum_{K \in \cB }  c_K h_K . $$
Observe that $ \|f \|_{\BMO} = 1 .$ Indeed, for $I_0 \in \cB $
write 
$$
\begin{aligned}
 \sum_{K \in \cB ,\, K \sbe I_0 } |K|  c_K^2  
&=\sum_{K \in \cB ,\, K \sbe I_0 } |K|- |G^*_1( K | \cB )| \\
&= |I_0| .
\end{aligned}
$$
Next observe that the right hand side of \eqref{2812so2} 
coincides with $ \|S_\s f\|_2^2 . $ We have 
$S_\s f = \sum_{K \in \cB }  c_K h_{\s( K)} , $ and 
\begin{equation} \label{2812so3}
\begin{aligned}
 \| \sum_{K \in \cB }  c_K h_{\s ( K)} \|^2_2 
              & = \sum_{K \in \cB } |\s(K)|  c_K^2 \\
              & = \sum_{K \in \cB }   | \sigma (K) |  \frac{|K|- |G^*_1( K | \cB )| }{|K|} . 
\end{aligned}
\end{equation}
The Haar support of $S_\s f$ is $ \s(\cB)$ which is contained in $ \s(\cH ). $
Hence
$$ \|S_\s f\|_2^2 \le |\s(\cH )^*| \|S_\s f\|_{\BMO}^2. $$
Combining \eqref{2812so2} and  \eqref{2812so3} with  the fact that $ \|f \|_{\BMO} = 1 $ we obtain
$$
\begin{aligned}
\int_0^1\mu_{\cH } (t) dt & \le   |\s(\cH )^*|\|S_\s f\|_{\BMO}^2 \\
                       & \le  |\s(\cH )^*|\|S_\s \|_{\BMO}^2. 
\end{aligned}
$$
This proves the right hand side estimate.

Next we turn to the left hand side estimate using $H^1- \BMO $ duality.
We prove that 
$$ \|T_{\t , 1 } \|_1^2 \le A_1C_1, $$
where  $A_1$ is the constant appearing in the atomic decomposition for $H^1 $ and  
$$ C_1 = \sup_{ \cH \sbe \t( \cD )} \frac{1}{ |\s(\cH)^*|}
 \int_0^1 
\mu _{\cH } ( t ) 
dt \quad \text{with}\quad \mu _{\cH } ( t ) = \sup _{J \in \cH } \frac{ |\s (J) |} {|J|} 1_{J }(t)
.
$$
Fix $ I  \in \cD $ and let $f : [0,1] \to \bR $ be a dyadic $H^1 $ atom so that 
\begin{equation} \label{2812so4}
 \supp S^2 (f) \sbe I \quad \text{and}\quad \|f \|_2 \le |I|^{-1/2}.
\end{equation}
The   square functions $S(  T_{\t, 1} (f) ) $ and $S(  T_{\t,2 }  (f))$
are related  by  pointwise factorization,
$$
S(  T_{\t, 1} (f) )  \le \mu _{\cH } ^{1/2} S(  T_{\t, 2} (f) )
$$
where $  \cH = \t(\{ J \in \cD : J \sbe I_0 \}  ) .$
Integrating and using the Cauchy Schwarz inequality gives,
\begin{equation} \label{2812so5} 
 \begin{aligned}
\|T_{\t, 1} (f)\|_{H^1}  &\le \left( \int_0^1\mu _{\cH } \right )^{1/2} \| T_{\t, 2} (f) \|_2\\
                            & \le C_1^{1/2}   |\s(\cH )^*|^{1/2} \| f \|_2. 
\end{aligned}
\end{equation}
Since $ \s(\cH ) ^* \sbe I$ and $f $ is a dyadic atom satisfying
\eqref{2812so4} we get from \eqref{2812so5} that 
$  \|T_{\t, 1} \|_{H^1}  \le  (A_1C_1)^{1/2}  $
and by duality,
$$ \| S_\s \|_{\BMO}^2 \le A C_1 . $$
\endproof
\paragraph{Remark.} Compare the inequalities of Proposition~\ref{proposition}
with \eqref{3feb1}. For
$$ C_1 = \sup_{ \cH \sbe \t( \cD )} \frac{1}{ |\s(\cH)^*|}
 \int_0^1 
\mu _{\cH } ( t ) 
dt \quad \text{with}\quad \mu _{\cH } ( t ) = \sup _{J \in \cH } \frac{ |\s (J) |} {|J|} 1_{J }(t)$$
we get
        $$C_1 \sim  \sup \frac{\y\s(\cC) \yy }{ \y\cC \yy} ,$$
where the supremum is extended over all $\cC \sbe \t ( \cD) .$

\begin{satz} \label{28jan091}Let $ 0<p \le 2. $ 
Then for any $0 < q \le  p,  $
$$ \|  T_{\t, q } \otimes \Id_X \| _{q} \le  A_q C_1^{1/q-1/p}
\|  T_{\t, p } \otimes \Id_X \| _{p} ,
$$
where 
 $$
C_1 = \sup_{ \cH \sbe \t( \cD )} \frac{1}{ |\s(\cH)^*|}
 \int_0^1 
\mu _{\cH } ( t ) 
dt  
$$
and $A_q>0$ is determined by the atomic decomposition  for  $H^q_X . $
\end{satz}
\proof
Let $ 0 < q \le p \le 2. $ Let $f$
be an  $X$ valued $H^q  $ atom,  so that 
$$ \supp \bS(f) \sbe I \quad\text{ and }\quad \bS(f) \le |I|^{-1/q} . $$
Comparing the defining equations for  $T_{\t, p } \otimes \Id_X (f)$
and $ T_{\t, q } \otimes \Id_X (f) , $
gives the pointwise estimate between square functions,
$$
\bS( T_{\t, q } \otimes \Id_X (f))(t)
\le 
\sup_{I \in \cL}  \left[\frac{ |I|}{ |\t (I)|}1_{\t(I)}(t)\right]^{1/q - 1/p}
\bS( T_{\t, p } \otimes \Id_X (f))(t). $$
Hence with $ \cH = \t^{-1} (\{ J \in \cD : J \sbe I \} ) $ we get the factorization 
$$
\bS( T_{\t, q } \otimes \Id_X (f))
\le \mu_{\cH} ^{1/q - 1/p}
\bS( T_{\t, p } \otimes \Id_X (f)). $$
Next raise the above estimate to the power $q$ and apply Hoelder's inequality 
with conjugate indices $ p/q$ and $ p/(p-q) . $ This gives
$$
\int \bS( T_{\t, q } \otimes \Id_X (f))^q dt 
\le \left(\int \mu_{\cH}(t) dt \right)  ^{1-q/p}
\left(\int \bS( T_{\t, p } \otimes \Id_X (f))^p dt\right)^{q/p}. $$
Taking $q -$th root yields 
$$ \|  T_{\t, q } \otimes \Id_X (f) \|_{H^q_X} \le 
\left(\int  \mu_{\cH}(t) dt \right)  ^{1/q-1/p}
\|  T_{\t, p} \otimes \Id_X (f) \|_{H^p_X} .$$
Since $\int \mu _{\cH } \le C_1 |I| $ and $\|f\|_{H^p_X} \le |I|^{1/p -1/q} , $
we get 
$$ \|  T_{\t, q } \otimes \Id_X (f) \|_{H^q_X} \le  C_1  ^{1/q-1/p}
\|  T_{\t, p} \otimes \Id_X \|_{p} 
.$$
The atomic decomposition theorem for $H^q_X$ implies now 
$$ \|  T_{\t, q } \otimes \Id_X \|_{q} \le A_q  C_1  ^{1/q-1/p}
\|  T_{\t, p} \otimes \Id_X \|_{p} 
.$$
\endproof 
In the case  when $ p>0 $ is strictly less that $2,$ the conclusion  of the previous theorem
can be turned into  concise extrapolation estimates as follows.
\begin{satz} \label{28jan0911}Let $ 0<q \le p <2. $ Then
$$ \|  T_{\t, q } \otimes \Id_X \| _{q}^{\frac{q}{2-q}}  \le 
A(p,q)\|  T_{\t, p } \otimes \Id_X \| _{p}^{\frac{p}{2-p}} .
$$
\end{satz}
\proof
Let  $ 0 < p < 2$ (and p strictly less than $2$). By \eqref{15jan081} 
we have the 
scalar valued extrapolation estimate \cite{gmp}
$$ \|T_{\t, 1}\|_1 \le A(p)\|T_{\t, p}\|_p^{\frac{p}{2-p}}.$$
Recall that $ C_1 \le (A\|T_{\t, 1}\|_1)^2.$ Hence the above inequality gives
$$  C_1^{\frac{p-q}{pq}} \le (A\|T_{\t, p}\|_p)^{\frac{2(p-q)}{p(2-p)}} .$$
Next note that $\frac{2(p-q)}{p(2-p)} + 1 = \frac{p(2-q)}{q(2-p)}. $
It remains to invoke  Theorem~\ref{28jan091} to obtain
$$
\begin{aligned}
 \|  T_{\t, q } \otimes \Id_X \| _{q} &\le A(p,q) C_1^{\frac{p-q}{pq}}
\|  T_{\t, p } \otimes \Id_X \| _{p}\\
&\le A(p,q) \|  T_{\t, p } \otimes I_X \| _{p}^{ \frac{p(2-q)}{q(2-p)}}.
\end{aligned}
$$
as claimed.
\endproof
\paragraph{Remark.}
Our proof  identifies the separate contribution  of
$\|  T_{\t, p } \|_p$ and  $\|  T_{\t, p } \otimes \Id_X \| _{p}$ 
to the upper bound for $\|  T_{\t, q } \otimes \Id_X \| _{q}. $
It gives
\begin{equation}\label{15jan084}
  \|  T_{\t, q } \otimes \Id_X \| _{q}\le A(p,q) \|  T_{\t, p } \|_p^{2\left(\frac{1}{q} - \frac{1}{p}\right)}
\|  T_{\t, p } \otimes \Id_X \| _{p}.
\end{equation}

The extrapolation estimates of Theorem \ref{28jan091} and Theorem 
\ref{28jan0911}  hold  for one fixed Banach space $X.$
Now we change the nature of our assumptions and demand   boundedness 
of the vector valued rearrangement operator for {\em each }
Banach space with the UMD property. While  this formulates a more restrictive  hypothesis
the resulting conclusion is also  stronger. 
The next theorem is a consequence of 
Proposition~\ref{8feb2} below and of the extrapolation theorems in  \cite{gmbull}.
We point out that in the hypothesis of Theorem~\ref{8feb1} the value $ p_0 = 2 $ is included.
\begin{satz} \label{8feb1}If there exists $1 < p_0 < \infty $  so that 
\begin{equation} \label{8feb3} 
\|  T_{\t, p_0 } \otimes \Id_E \| _{L^{p_0}_E} < \infty ,\quad\text{for each UMD space E.}
\end{equation}
then  $\|T_{\tau,1}^{-1}\|_1 <\infty  $ and 
for any  $1 < p < \infty, $
$$ \|  T_{\t, p } \otimes \Id_E \| _{L^{p}_E} \|  T_{\t, p }^{-1} \otimes \Id_E \| _{L^{p}_E} < \infty, $$
for each UMD space $E.$
\end{satz}

Consider a rearrangement operator 
(of the $H^1$ normalized Haar system) that is unbounded on $H^1 . $ Then
as shown in \cite{pfxm}, for any $n,$ 
there exist vectors $ x_ 1 , \dots , x_n $ 
in $H^1 $ that are equivalent to the unit vector basis of $\ell^2_n ,$
 so that their images are equivalent to the unit vector basis of $\ell^1_n .$
Guided by the reasoning of \cite[Example 3.2]{gmpre}
we next give the vector valued interpretation 
of Proposition 2 and Lemma 3 +4 in \cite{pfxm}. 
Thus we show that \eqref{8feb3} implies $\|T_{\tau,1}^{-1}\|_1 <\infty . $
\begin{prop}\label{8feb2}
For any Banach space $Y$ and $ 1 < p \le 2 $ we have
$$\|T_{\tau,1}^{-1}\|_1 = \infty \quad\text{implies}\quad \text{Type} _{p} (Y)\le C \|  T_{\t, p } \otimes \Id_Y \| _{L^{p_0}_Y} .$$
Consequently, if 
$$\|  T_{\t, p } \otimes \Id_E \| _{L^{p}_E} < \infty ,\quad\text{for each UMD space E}, $$
then 
$$\|T_{\tau,1}^{-1}\|_1 < \infty . $$
\end{prop}
\proof
Suppose that $ \| T_{\tau, 1}^{-1} \|_1 = \infty . $
By \cite{pfxm} (see also Proposition 3.3.2 and Theorem 3.3.5 in \cite{m}) 
for each $n\in \bN $ there exists a collection of dyadic intervals 
$\cC $ so that 
$$\y \cC \yy \le 4\quad\text{and}\quad \y \t( \cC ) \yy \ge n^2 . $$
By the Carleson-Garnett
condensation lemma, (see \cite[Lemma 3.1.4]{m}), 
there exists $ K \in \t( \cC ) $ so that 
\begin{equation}\label{8feb7}
 |G_ n (K | \t ( \cC ))^* | \ge ( 1 - \frac1n)|K|,
\end{equation}
where $G_ n (K | \t ( \cC ))$ denotes the $n -$th generation of $\t(\cC) $ that is underneath $K. $
By rescaling we may assume that $K = [ 0 , 1 ] . $ Now  for $ i \le n ,$ put
$$ \cF_i = G_ i (K | \t ( \cC ))\quad\text{and} \quad \cE_i = \t^{-1} G_ i (K | \t ( \cC ))$$
to define 
$$ \rho _ i = \sum_{ \t(J) \in  \cF_ i} h_{\t(J)}  \quad \text{and} 
\quad s_i = \sum_{ J \in \cE_i }  | \t(J)|^{1/p}\frac{  h_J }{ |J|^{1/p} } h_J .$$
Note that 
$$
 \bigcup_{i = 1 } ^n \cE_i \sbe \cC \quad\text{and} \quad T_{\tau, p}s_i = \rho _i  \quad i \le n . $$
Since  $\y \cC \yy \le 4$ we have with \cite[Lemma 3]{gmbull} that for $a_i \in Y$ 
\begin{equation}\label{8feb30}
\begin{aligned} 
\| \sum_{i = 1 }^ n a_i s_i \|_{L^p_Y} & \le C \left( \sum_{i = 1 }^ n \|a_i \|_Y^p \sum_{J \in \cE_i} |\t ( J )|  
\right)^{1/p}\\
& \le C\left( \sum_{i = 1 }^ n \|a_i \|_Y^p
\right)^{1/p}.
 \end{aligned}
\end{equation}

Now let $ \{ r_i \} _{i = 1 } ^ n$ be the first $n $ of the independent $ \{ \pm 1 \} $ valued 
Rademacher functions.  It follows from \eqref{8feb7} that   for $a_i \in Y$ 
\begin{equation}\label{8feb33}
\| \sum_{i = 1 }^ n a_i r_i \|_{L^p_Y}  \le \| \sum_{i = 1 }^ n a_i \rho_i \|_{L^p_Y} + 
\frac{1}{n}\left( \sum_{i = 1 }^ n \|a_i \|_Y\right) .
\end{equation}
Hence for $ n \in \bN $ large enough, by \eqref{8feb30} and \eqref{8feb33}, 
we get with $ \rho _i = T_{\tau, p}s_i$ 
\begin{equation}\label{8feb32}
\begin{aligned}
\| \sum_{i = 1 }^ n a_i r_i \|_{L^p_Y} &  \le   \| \sum_{i = 1 }^ n a_i T_{\tau, p}s_i \|_{L^p_Y} + 
\frac{1}{n}\left( \sum_{i = 1 }^ n \|a_i \|_Y\right)  \\
& \le C \|  T_{\t, p } \otimes \Id_Y \| _{L^{p}_Y}\| \sum_{i = 1 }^ n a_i s_i \|_{L^p_Y}+
\frac{1}{n}\left( \sum_{i = 1 }^ n \|a_i \|_Y\right) \\
& \le C\|  T_{\t, p } \otimes \Id_Y \| _{L^{p}_Y} \left( \sum_{i = 1 }^ n \|a_i \|_Y^p
\right)^{1/p} . 
\end{aligned}
\end{equation}
Since $ i \le n $ and  $ a_i \in Y $ were chosen arbitrary \eqref{8feb32} 
implies
$$\text{Type} _{p} (Y)\le C  \|  T_{\t, p } \otimes \Id_Y \| _{L^{p}_Y} .$$ 
To see the moreover part 
of Proposition~\ref{8feb2} just test the above estimate 
with $ Y = L^ r , $ and $ 1 <r < p . $
\endproof

\paragraph{Proof  of Theorem~\ref{8feb1}.} It suffices to consider $ 1 < p_0 \le 2.$ Consider first the case $p_0 = 2 . $
 Since UMD spaces are reflexive, and the UMD property is a self dual isomorphic invariant we get 
from \eqref{8feb3} by dualization that  
\begin{equation}\label{8feb20}
\|  T_{\t, 2 } \otimes \Id_E \| _{L^{2}_E}\|  T_{\t, 2 }^{-1} \otimes \Id_E \| _{L^{2}_E} < \infty ,
\end{equation}
for any UMD space $E .$
Hence by Proposition~\ref{8feb2}
\begin{equation}\label{8feb21}
\|T_{\tau,1}^{-1}\|_1\|T_{\tau,1}\|_1 < \infty .
\end{equation}
By \cite[Corollary 5.6.]{gmpre},  we get from  \eqref{8feb20} and \eqref{8feb21}
that for any  $1 < p \le 2 $ and any UMD space $E$
$$ \|  T_{\t, p } \otimes \Id_E \| _{L^{p}_E} \|  T_{\t, p }^{-1} \otimes \Id_E \| _{L^{p}_E} < \infty. $$
By duality this gives  the conclusion of Theorem~\ref{8feb1} in the case  $p_0 = 2. $

Next we turn to  $p_0 <  2. $
By scalar valued extrapolation we get then $\|T_{\tau,1}\|_1 < \infty , $ and  Proposition~\ref{8feb2}
yields $\|T_{\tau,1}^{-1}\|_1 < \infty . $ Hence as in the previous case
\begin{equation}\label{8feb24}
\|T_{\tau,1}^{-1}\|_1\|T_{\tau,1}\|_1 < \infty . \end{equation}
Moreover by reflexivity, duality, and the fact that UMD is a self dual 
isomorphic invariant, 
 we obtain  with $1/p_0 + 1/q_0 = 1 ,$ that 
\begin{equation}\label{8feb25}
\|  T_{\t, q_0 }^{-1} \otimes \Id_E \| _{L^{q_0}_E}< \infty
\quad\text{for each UMD space E.
}
 \end{equation}
By \cite[Corollary 5.6.]{gmpre} it follows from   \eqref{8feb24} and \eqref{8feb25} that 
for each $ q \le q_0 $, the operator   $T_{\t, q }^{-1} \otimes \Id_E$ is 
bounded on 
${L^{q}_E}.$  Since $ p_0 < 2 $ and $ q_0 >  2 , $ this holds in particular for $ q = p_0 . $
In summary we have 
$$ \|  T_{\t, p_0 } \otimes \Id_E \| _{L^{p_0}_E} \|  T_{\t, p_0 }^{-1} \otimes \Id_E \| _{L^{p_0}_E} < \infty . $$
  By \cite[Corollary 5.8.]{gmpre} we obtain the conclusion of Theorem~\ref{8feb1}
in the case $p_0 < 2.$
\endproof

\nocite{gmpre}
\nocite{maureyfactorisation} 
\nocite{maureyvancouver}
\nocite{bu1}
\nocite{blasco}
\nocite{gmp}
\nocite{pfxm}
\nocite{gmbull}
\bibliographystyle{abbrv}
\bibliography{tot}   
\paragraph{Address:}   
Institut f\"ur Analysis \\ J. Kepler Universit\"at\\A4040 Linz Austria \\
pfxm@bayou.uni-linz.ac.at
\end{document}